\documentclass[a4paper,10pt]{article}
\usepackage[T1]{fontenc}
\usepackage{amsmath}
\usepackage{amssymb}
\usepackage{graphicx}
\usepackage{color}
\usepackage{bbold}

\usepackage{hyperref}

\renewcommand{\P}{{\rm P}}
\newcommand{\E}{{\rm E}}

\newcommand{\indic}{\mathbb{1}}

\newcommand{\tr}{^{\textnormal {\tiny  T}}}

\newcommand{\ud}{\mathrm{d}}  

\newcommand{\s}{\mathcal{S}}
\newcommand{\cald}{\mathcal{D}}

\newcommand{\vect}[1]{\boldsymbol #1}

\newcommand{\vxi}{\vect \xi}
\newcommand{\vpi}{\vect \pi}

\newcommand{\vomega}{\vect \omega}

\newcommand{\vone}{\vect 1}

\newcommand{\vligne}[1]{\begin{bmatrix} #1 \end{bmatrix}}

\newtheorem{defn}{Definition}[section]

\newtheorem{lemma}[defn]{Lemma}

\newtheorem{theorem}[defn]{Theorem}

\newtheorem{cor}[defn]{Corollary}

\newtheorem{remark}[defn]{Remark}

\newtheorem{example}[defn]{Example}

\newcommand{\qed}{\hfill $\square$}
\newenvironment{proof}{
     \noindent {\bf Proof }}{\qed
     \vspace{0.25\baselineskip}}
\newcommand{\debproof}{\begin{proof}}
\newcommand{\finproof}{\end{proof}}

\newcommand{\passe}[1]{\vspace{#1\baselineskip}}

\definecolor{darkmagenta}{rgb}{0.5,0,0.5}
\definecolor{darkgreen}{rgb}{0,0.6,0}
\definecolor{darkblue}{rgb}{0,0,0.6}
\definecolor{darkred}{rgb}{0.8,0,0}
\definecolor{mellow}{rgb}{.847, 0.72, 0.525}

\usepackage{pifont}

\usepackage[normalem]{ulem}

\begin{document}

\title{Why is Kemeny's constant a constant?}

\author{
Dario Bini\thanks{University of Pisa, Dipartimento di Matematica,
  56127 Pisa, \texttt{bini@dm.unipi.it, meini@dm.unipi.it}}
\and
Jeffrey J. Hunter\thanks{Auckland University of Technology, Department of Mathematical Sciences, 1142 Auckland, \texttt{jeffrey.hunter@aut.ac.nz}}
\and
Guy Latouche\thanks{Universit\'e Libre de Bruxelles,  D\'epartement
  d'informatique, 1050 Bruxelles, \texttt{latouche@ulb.ac.be}}
\and
Beatrice Meini$^*$
\and
Peter Taylor\thanks{University of Melbourne, School of Mathematics and Statistics, Vic 3010, \texttt{taylorpg@unimelb.edu.au}}
}


\maketitle

\begin{abstract}
In their 1960 book on finite Markov chains, Kemeny and Snell established that a certain sum is invariant.  The value of this sum
has become known as {\it Kemeny's constant}.  Various proofs have been given over time, some more technical than others.  We give here a very simple physical justification, which extends without a hitch to continuous-time Markov chains on a finite state space. 

For Markov chains with denumerably infinite state space, the constant may be infinite and even if it is finite, there is no guarantee that the physical argument will hold. We show that the physical interpretation does go through for the special case of a birth-and-death process with a finite value of Kemeny's constant.

\passe{0.5}

\passe{0.5}
\end{abstract}

\section{Introduction}
\label{s:introduction}

Consider a discrete-time, irreducible and aperiodic Markov chain $\{X_t: t= 0, 1,
\ldots\}$ on a finite state
space $\s$, with transition matrix $P$ and stationary probability
vector $\vpi$ such that $\vpi\tr P = \vpi\tr$ and $\vpi\tr \vone = 1$.
For $i \in \s$, define the {\em first passage times} 
\begin{equation}
   \label{e:return}
T_i   = \inf\{ t \geq 1: X_t=i\}.
\end{equation}
Denoting by $\E_i[\cdot]$ the conditional expectation given that $X_0=i$, Kemeny and Snell~\cite[Theorem 4.4.10]{ks60} proved that
\begin{equation}
   \label{e:K}
\sum_{j \in \s} \pi_j \E_i[T_j] = K,
\end{equation}
independently of the initial state $i$. The value $K$ is known as {\it
  Kemeny's constant}.

A prize was offered to the first person to give an intuitively
plausible reason for the sum in (\ref{e:K}) to be independent of $i$
(Grinstead and Snell~\cite[Page 469]{gs97}). The prize was won by
Doyle~\cite{doyle09} with an argument given in the next section.  We
prove in Theorem~\ref{t:kemenyfinite} that (\ref{e:K}) results from
the obvious fact that a discrete-time Markov chain takes $n$ steps
during an interval of time of length $n$, independently of the initial
state $i$. We move on to extend the argument to finite-state {\em
  continuous-time} Markov chains, see (\ref{e:indepc}). In Section
\ref{s:finite}, we also discuss an important connection between $K$
and the deviation matrix of the Markov chain.

In Section~\ref{s:infinite}, we consider Markov chains with a denumerably infinite state space $\s$. Here, the situation becomes more complex because the sum in (\ref{e:K}) might not converge. We show that it is independent of $i$ in the sense that it is infinite for all $i$ or a constant independent of $i$. 

In Section~\ref{s:bandd} we restrict our discussion to positive recurrent birth-and-death  processes.   We show that $K$ is infinite in discrete-time, and in continuous-time it is finite if transitions from state $i$ occur sufficiently fast as $i$ approaches infinity. Furthermore, our physical explanation holds for birth-and-death processes if $K < \infty$.

\section{A simple algebraic proof}
\label{s:simple}

The simplest proof goes as follows: define $\omega_i = \sum_{j
  \in \s} \pi_j \E_i[T_j]$ and $\vomega= \vligne{\omega_i}_{i \in \s}$, condition on $X_1$ and write
\begin{align*}
\omega_i 
 & = 1 + \sum_{j \in \s} \pi_j \sum_{k \in \s, k \not= j}
   P_{ik}\E_k[T_j] \\
 & = 1 + \sum_{j \in \s} \pi_j \sum_{k \in \s}    P_{ik}\E_k[T_j]  -
   \sum_{j \in \s} P_{ij},  \qquad \mbox{using $\pi_j = 1/ \E_j[T_j]$,}\\
  & = \sum_{j \in \s} \pi_j \sum_{k \in \s}    P_{ik}\E_k[T_j] \\
  & = \sum_{k \in \s} P_{ik} \,\omega_k
\label{eq:doyle}
\end{align*}
so that $\vomega = P \vomega$ (see, for example, Hunter \cite{hunt14}).
Doyle~\cite{doyle09} argued from the maximum principle that all
components of $\vomega$ must be equal.  Alternatively, one may conclude from the Perron-Frobenius Theorem that $\vomega$ must be proportional to the eigenvector $\vone$ of $P$.

Instead of the passage times $T_j$, we shall use the {\em first hitting
  times}  $\{\theta_i : i \in \s\}$ with
\begin{equation}
   \label{e:passage}
\theta_i  = \inf\{ t \geq 0: X_t=i\}.
\end{equation}
The only difference is that $\theta_i = 0 < T_i$ if $X_0 = i$, 
otherwise $\theta_i = T_i \geq 1$.  Using $\theta_j$ instead of $T_j$,
we obtain another version of Kemeny's constant:
\begin{equation}
   \label{e:kprime}
\sum_{j \in \s} \pi_j \E_i[\theta_j] = K'
\end{equation}
where $K'=K-1$.  We prefer to work with this version of Kemeny's
constant because the equality (\ref{e:kprime}) holds in
continuous-time as well; furthermore, using $\E_i[\theta_j]$ helps us
establish a direct connection with the deviation matrix $\cald$ of the
Markov chain.  We shall discuss this in the next section.

\section{The case when $\s$ is finite}
  \label{s:finite}

Our physical justification is based on the following argument.  We
start from
\begin{equation}
   \label{e:start}
\sum_{j \in \s} \pi_j \E_i[\theta_j]  = \sum_{j \in \s}  \frac{\E_i[\theta_j] }{\E_j[T_j]}
\end{equation}
which we transform to
\begin{equation}
   \label{e:diffNs}
\sum_{j \in \s} \pi_j \E_i[\theta_j]  =  \sum_{j \in \s} \lim_{n \rightarrow \infty}  (\E_j[N_j(n)] - \E_i[N_j(n)]),
\end{equation}
where 
\[
N_j(n) = \sum_{0 \leq t \leq n} \indic\{X_t=j\}.
\]
is the total number of visits to $j$ during the interval of time $[0,n]$.

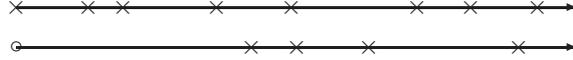
\begin{figure}[ht!]
\begin{center}


\begin{picture}(210,25)(0,0)
\put(0,20){\vector(1,0){210}}
  \put(0,20){\makebox(0,0){$\times$}}
  \put(27,20){\makebox(0,0){$\times$}}
  \put(40,20){\makebox(0,0){$\times$}}
  \put(75,20){\makebox(0,0){$\times$}}
  \put(103,20){\makebox(0,0){$\times$}}
  \put(150,20){\makebox(0,0){$\times$}}
  \put(170,20){\makebox(0,0){$\times$}}
  \put(195,20){\makebox(0,0){$\times$}}
\put(0,5){\vector(1,0){210}}
  \put(0,5){\makebox(0,0){$\circ$}}
  \put(88,5){\makebox(0,0){$\times$}}
  \put(105,5){\makebox(0,0){$\times$}}
  \put(132,5){\makebox(0,0){$\times$}}
  \put(188,5){\makebox(0,0){$\times$}}
\end{picture}

\caption{Visits to $j$ starting in $j$ (above) and $i$ (below).}
\label{fig:samplepaths}
\end{center}
\end{figure}

The formal justification for the transition from (\ref{e:start})
to (\ref{e:diffNs}) is given in Lemma~\ref{t:newview} below, but we give
a heuristic argument first, explained with the help of Figure \ref{fig:samplepaths}.
The upper line is a representation of a trajectory of the renewal
process $\{ \theta_j^{(k)} : k \geq 0 \} $ of successive visits to
$j$, starting from $X_0=j$; the $\theta_j^{(k)}$s are marked with a cross
$\times$.  The lower line represents a trajectory of the delayed
renewal process of visits to $j$, starting from $X_0=i\not = j$.  

Now, the $j$th term in the right-hand side of (\ref{e:diffNs}) is the
expected difference between the total number of events in the two processes.
We observe a smaller expected number of visits to $j$ if the process starts from $i\not= j$ because of the
initial delay. The expected length of this delay is $\E_i[\theta_j]$ and
$\E_j[T_j]$ is the expected length of intervals between visits to $j$. The ratio $\E_i[\theta_j] /\E_j[T_j]$ is the expected number of
visits that are missed over the whole history of the process by
starting from $i$ instead of $j$.  The formal argument is given now.

\begin{lemma}
   \label{t:newview}
For all $i$ and $j$.
\begin{equation}
\frac{\E_i[\theta_j] }{\E_j[T_j]} = \lim_{n \rightarrow \infty}
(\E_j[N_j(n)] - \E_i[N_j(n)]).
\label{e:lemma1}
\end{equation}

\end{lemma}
\begin{proof}
  The statement is obvious if $i=j$ for then $\E_j[\theta_j]=0$ by
  the definition of $\theta_j$.  We assume now that $i$ and $j$ are
  different, arbitrary but fixed, and to simplify the notation we
  define $\widetilde N_i(n)= \E_i[N_j(n)]$ and $f_i(t)= \P_i[\theta_j=t]$,
  with $f_i(0)=0$.  We have
\[
\widetilde N_j(n) = \sum_{0 \leq \nu \leq n} \P_j[X_\nu = j].
\] 
Furthermore, conditioning on the first visit to state $j$, we can write, for $n \geq 0$,
\begin{align*}
\widetilde N_i(n) & = \sum_{0 \leq t \leq n} f_i(t)  \widetilde
                    N_j(n-t) \\
 & = \sum_{0 \leq t \leq n} f_i(t) \sum_{0 \leq \nu \leq n-t}
   \P_j[X_\nu = j] \\
 & = \sum_{0 \leq \nu \leq n} \P_j[X_\nu = j] \sum_{0 \leq t \leq n-
   \nu} f_i(t)
\end{align*}
 and so
\begin{align*}
 \widetilde N_j(n)- \widetilde N_i(n) & = \sum_{0 \leq \nu \leq n}
                                        \P_j[X_\nu = j]  \P_i[\theta_j > n-\nu].
\end{align*}
Finally, 
\begin{align*}
\lim_{n \rightarrow \infty }    (\widetilde N_j(n)- \widetilde N_i(n) )
 & = \lim_{n
   \rightarrow \infty }   \sum_{0 \leq \nu \leq n} \P_j[X_{n-\nu}=j] \P_i[\theta_j
                                                       > \nu]   \\
  & = \frac{1 }{\E_j[T_j]}   \sum_{\nu \geq 0}   \P_i[\theta_j> \nu]
\intertext{by the key renewal theorem (Resnick~\cite[Section~3.8]{resni92}),}  
  & = \pi_j \E_i[\theta_j].
\end{align*}
This completes the proof. 
\end{proof}

Lemma \ref{t:newview} leads immediately to a understanding of the reason why the left hand side of (\ref{e:K}) is independent of $i$.

\begin{theorem}
\label{t:kemenyfinite}
For an irreducible and aperiodic discrete-time Markov chain with finite state space $\s$ and stationary distribution $\vpi$, $\sum_{j \in \s} \pi_j \E_i[T_j]$ is independent of $i$.
\end{theorem}

\begin{proof}
Since $\s$ is finite, we can interchange the limit and sum in Equation (\ref{e:diffNs}) to get
\begin{align}
  \nonumber
\sum_{j \in \s} \pi_j \E_i[\theta_j]   & = 
   \sum_{j \in \s} \lim_{n \rightarrow \infty} (\E_j[N_j(n)] - \E_i[N_j(n)]) \\
  \nonumber & = \lim_{n \rightarrow \infty}
   \left(\sum_{j \in \s} \E_j[N_j(n)] - \sum_{j \in \s} \E_i[N_j(n)]\right)  \\
   \label{e:indept}
 & = \lim_{n \rightarrow \infty} \left(\sum_{j \in \s} \E_j[N_j(n)] - (n+1)\right)
 \end{align}
independently of $i$.
\end{proof}

We see that the change of perspective from first hitting times to numbers of visits brings a different physical interpretation from those previously suggested. Indeed, as mentioned above, the second term under the limit in (\ref{e:indept}) just counts the number of steps that the Markov chain has taken, independently of its starting state.

The deviation matrix lends itself beautifully to such a change of point of view.  It is defined
as
\begin{equation}
   \label{e:deviation}
\cald = \sum_{n \geq 0} (P^n - \vone \cdot \vpi\tr)
\end{equation}
if the series converges.  By Syski~\cite[Proposition 3.2]{syski78}, and Coolen-Schrijner and van Doorn~\cite[Theorem 4.1]{cv02}, the series
converges if and only if $\E_{\vpi}[\theta_j] < \infty$ for some state
$j \in \s$, and then it is finite for every $j$,
where $\E_\pi[\cdot]$ denotes
the conditional expectation, given that $X_0$ has the distribution
$\vpi$.  An equivalent
condition is that $\E_j[\theta_j^2] < \infty$ for some $j$. When $|\s| < \infty$, the series always converges and $\cald = (I-P)^\#$, the group inverse of $I-P$ (see Campbell and Meyer~\cite{cm91}).

Obviously, 
\begin{align}
  \nonumber
\cald_{ij} & = \lim_{n \rightarrow \infty} \sum_{t=0}^n ([P^t]_{ij} - \pi_j) \\
  \nonumber & = \lim_{n \rightarrow \infty} \sum_{t=0}^n (E_i[I(X_t=j)] - E_\pi[I(X_t=j)]) \\
& = \lim_{n \rightarrow \infty} (\E_i[N_j(n)] - \E_\pi[N_j(n)]).
 \end{align}
In addition, $\cald_{jj} = \pi_j \E_\pi[\theta_j]$ by Syski~\cite[Proposition 3.3]{syski78} with the discrete-time analogue of the argument in Coolen-Schrijner and
van Doorn~\cite[Equation (5.7)]{cv02}, and so
\begin{align}
  \nonumber
\sum_{j \in \s} \cald_{jj} & = \sum_{j \in \s} \pi_j \E_\pi[\theta_j] \\
  \nonumber  & = \sum_{j \in \s} \pi_j \sum_{i \in \s} \pi_i \E_i[\theta_j] \\
  \nonumber & = \sum_{i \in \s} \pi_i \sum_{j \in \s} \pi_j \E_i[\theta_j] \\
& = \sum_{i \in \s} \pi_i K' = K'. 
 \label{e:kdev}
\end{align}
Equation (\ref{e:kdev}) provides a convenient representation for Kemeny's constant in terms of the trace of the deviation matrix.

The point of view that we have taken in this section extends to finite-state continuous-time Markov chains. For such a chain with irreducible generator
$Q$, define the first hitting time and the first passage time as
\begin{align*}
\theta_i & = \inf\{ t \geq 0: X_t=i\} \\
T_i  & = \inf\{ t \geq J_1: X_t=i\},
\end{align*}
where $J_1$ is the first jump time of the Markov chain;
if $X_0 = i$, then $\theta_i = 0 < T_i$, otherwise $\theta_i = T_i
> 0$.
Lemma \ref{t:newview} becomes

\begin{lemma}
   \label{t:continuous}
For a continuous-time Markov chain, 
\[
\pi_j \E_i[\theta_j] = 
  \lim_{t \rightarrow \infty} (\E_j[M_j(t)] - \E_i[M_j(t)])
\]
where $M_j(t) = \int_0^t \indic\{X(u)=j\} \, \ud u$ is the total time
spent in $j$ until time $t$.
\end{lemma}
\begin{proof}
We follow the same steps as in Lemma~\ref{t:newview} with the only
difference that here
\[
\E_j[M_j(t)] = \int_0^t \P_j[X(u)=j] \, \ud u
\] 
and 
\[
\E_i[M_j(t)] = \int_0^t  \ud G_i(v) \E_j[M_j(t-v)]
\]
where $G_i(t) = \P_i[\theta_j \leq t]$.
\end{proof}

From this, we obtain 
\begin{equation}
\sum_{j \in \s} \pi_j \E_i[\theta_j] = \lim_{t \rightarrow \infty} 
  (\sum_{j \in \s}\E_j[M_j(t)] - t) 
\label{e:indepc}
\end{equation}
independently of $i$. By an argument similar to that which led to (\ref{e:kdev}), we also have
\begin{align}
\sum_{j \in \s} \pi_j \E_i[\theta_j]   & = \sum_{j \in \s} \pi_j \E_{\vpi} [\theta_j]
  \nonumber \\
  & = \sum_{j \in \s}  \cald_{jj}
  \label{e:Kcontinuous}
\end{align}
with the continuous-time deviation matrix defined by $\cald = \int_0^\infty (e^{Qt} - \vone \cdot
\vpi\tr) \, \ud t$.

\section{The case when $\s$ is infinite}
  \label{s:infinite}

If $\s$ is denumerably infinite, it is not easy to see in general how the exchange of limit and sum inherent in the step between the first and second equations of (\ref{e:indept}) can be justified. However, $\sum_{j \in \s} \pi_j \E_i[\theta_j]$ is still independent of $i$ in the sense that it is either finite and constant with respect to $i$, or infinite for all $i$.

\begin{theorem}
   \label{t:infinite}
For an irreducible,  positive-recurrent discrete or continuous-time Markov chain $\{X_t\}$ with a countably-infinite state space $\s$, either 
\begin{enumerate}
\item $\sum_{j \in \s} \pi_j \E_i[\theta_j]$ is equal to a finite constant that is independent of $i$, or
\item $\sum_{j \in \s} \pi_j \E_i[\theta_j]$ is infinite for all $i\in\s$.
\end{enumerate}
\end{theorem}

\begin{proof}
The argument presented in Section 2 due to \cite{doyle09, hunt14} goes through even when the state space is infinite. Writing the expressions in terms of $E_i[\theta_j]$ rather than $E_i[T_j]$, for a discrete-time Markov chain with transition matrix $P$ we have
\begin{eqnarray}
\xi_i & \equiv & \sum_{j \in \s} \pi_jE_i[\theta_j] \nonumber\\
& = & \sum_{j \not = i} \pi_j \left[1 + \sum_{k \in \s} P_{ik} E_k[\theta_j]\right]\nonumber\\
& = & 1 - \pi_i  + \sum_{k \in \s} P_{ik} \sum_{j \not = i} \pi_j
      E_k[\theta_j] \nonumber\\	
& = & 1 - \pi_i  + \sum_{k \in \s} P_{ik} (\sum_{j \in \s} \pi_j E_k[\theta_j]
 - \pi_i   E_{k}[\theta_i]) \label{e:a} \\  
 & = & 1 - \pi_i  + \sum_{k \in \s} P_{ik}  \sum_{j \in \s} \pi_j E_k[\theta_j]    -
       \pi_i\sum_{k\in \s} P_{ik}  E_{k}[\theta_i]  \label{e:b}
\end{eqnarray}
where the series in (\ref{e:a}) and the first series in (\ref{e:b}) both
converge or both diverge, and the second series in (\ref{e:b})
converges to $E_{i}[T_i] -1 < \infty$ by assumption.
Thus, we may write
\begin{eqnarray}
\xi_i
& = & 1 - \pi_i  + \sum_{k \in \s} P_{ik} \xi_k -
      \pi_i\left[E_{i}[T_i]-1\right]\nonumber\\  
& = & \sum_{k \in \s} P_{ik} \xi_k.
\label{eq:doyle2}
\end{eqnarray}																
Since we have assumed that $\{X_t\}$ is recurrent, it follows from
Theorem 5.4$^D$ of Seneta \cite{sene81} that, if $\vxi$ is entrywise
finite, then it must be a multiple of $\vone$. On the other hand, if
$\xi_k$ is infinite for some $k$, then (\ref{eq:doyle2}) implies that $\xi_i$ must
be infinite for any $i$ such that $P_{ik}>0$. It follows by
irreducibility that  $\xi_i$ must be infinite for all $i \in \s$.

For a continuous-time Markov chain with transition matrix $Q$, similar reasoning holds with $P$ the transition matrix of the jump chain with entries $P_{ij} = q_{ij} I[i\not = j]/q_i$.
\end{proof}

\begin{remark}  \rm
\label{r:discrete}

In Section \ref{s:bandd}, we shall show that, for a discrete-time birth and death process with state space $\{0,1,\ldots\}$,  $K' = \sum_{j \in \s} \pi_j \E_i[\theta_j]$ is infinite for all $i$. 

For a general $m$-state discrete-time Markov chain, Hunter \cite[Theorem 4.2]{hunt06} used a spectral argument to show that $K\geq (m+1)/2$ which implies that $K'\geq (m-1)/2$. We do not see how to extend this argument to show that $K'$ is infinite for a general infinite-state discrete-time Markov chain, but we do not know of an example of such a chain with finite $K'$. We conjecture that $K'$ is infinite for all infinite-state, discrete-time Markov chains.

On the other hand, it is possible for $K'$ to be finite for an infinite-state continuous-time Markov chain. We shall present some examples in Section \ref{s:bandd}.

\end{remark}

\begin{remark}  \rm
\label{r:deviation}

In Section \ref{s:finite}, we showed that, when $\s$ is finite, $K'$ is equal to the trace of the deviation matrix. The argument in Appendix A shows that this is the case when $\s$ is infinite and $\cald$ exists. Specifically $K'$ is finite and equal to $\sum_{j \in \s} \cald_{jj}$ if this sum is finite, and infinite when $\sum_{j \in \s}  \cald_{jj}$ is infinite. 

When $\cald$ does not exist,  $\E_{\vpi}[\theta_j]$ must be infinite for all $j \in \s$, by Syski~\cite[Proposition 3.2]{syski78}, and Coolen-Schrijner and van Doorn~\cite[Theorem 4.1]{cv02}. Then
\begin{align}
\sum_{i \in \s} \pi_i \sum_{j \in \s} \pi_j \E_{i}[\theta_j] & = \sum_{j \in \s} \pi_j \sum_{i \in \s} \pi_i \E_{i}[\theta_j]\nonumber\\ 
& = \sum_{j \in \s} \pi_j \E_{\vpi}[\theta_j] \nonumber\\
& = \infty.  \nonumber
\end{align}
Since $\sum_{j \in \s} \pi_j \E_{i}[\theta_j]$ must be independent of
$i$ if it is finite, the last equality can occur only if $K'= \sum_{j \in \s} \pi_j \E_{i}[\theta_j] = \infty$.
\end{remark}
 
\section{Birth-and-death processes}
 \label{s:bandd}

Let us assume now that $\{X_t\}$ is a birth-and-death process on the
infinite state space $\{0, 1, \ldots\}$.  Choosing $X_0=0$ without
loss of generality, we have
\begin{equation}
   \label{e:sum}
K' = \sum_{j \geq 0} \pi_j \E_0[\theta_j].
\end{equation}
We shall examine continuous and discrete-time processes simultaneously.  In continuous-time, we denote by $\lambda_n$ and $\mu_n$ the transition rates from $n$ to $n+1$ and from $n$ to $n-1$, respectively; in discrete-time these are the one-step transition probabilities. We assume that $\lambda_n >0$ for all $n \geq 0$, $\mu_n >0$ for all $n \geq 1$, so
that the process is irreducible.    We further assume that {the birth-and-death process} is positive recurrent, so that $B = \sum_{n \geq 0} \beta_n$ is finite, where
\begin{equation}
   \label{e:betan}
\beta_0 = 1, \qquad \beta_n = \frac{\lambda_0 \lambda_1 \cdots
  \lambda_{n-1}}{\mu_1 \mu_2 \cdots \mu_n} \quad \mbox{for $n \geq 1$,}
\end{equation}
and the stationary distribution is given by $\pi_n = B^{-1} \beta_n$, {see
\cite[Equation~(6.3)]{cv02}, which is still valid in the discrete-time case.}

\begin{theorem}
   \label{t:bdinfinite}
For an irreducible, positive recurrent, birth-and-death process on
$\{0, 1, \ldots\}$, the constant $K'$ is finite if and only if 
\begin{equation}
   \label{e:G}
\Theta = \sum_{k \geq 0} (\lambda_k \pi_k)^{-1} \sum_{j \geq k+1} \pi_j < \infty.
\end{equation}
In that case, $K' = \Theta - \E_\pi [\theta_0]$, with
\begin{equation}
   \label{e:pitheta0}
\E_\pi [\theta_0] = \sum_{k \geq 0} (\lambda_k \pi_k)^{-1} (\sum_{j
  \geq k+1} \pi_j)^2 < \Theta.
\end{equation}
\end{theorem}
\begin{proof}
We start from (\ref{e:sum}) and write
\begin{align}
   \nonumber
K' & = \sum_{j \geq 1} \pi_j \sum_{0 \leq k \leq j-1} (\lambda_k \pi_k)^{-1}
\sum_{0 \leq \ell \leq k} \pi_\ell   \\
   \nonumber 
  & = \sum_{k \geq 0} (\lambda_k \pi_k)^{-1}
(\sum_{j \geq k+1} \pi_j )
(1 - \sum_{ \ell \geq k+1} \pi_\ell)
 \\
   \label{e:one}
  & = \sum_{k \geq 0} (\lambda_k \pi_k)^{-1}
\sum_{j \geq k+1} \pi_j      -
\sum_{k \geq 0} (\lambda_k \pi_k)^{-1} (\sum_{j \geq k+1} \pi_j )^2
\end{align}
if both series converge, with the first equation following from \cite[Equation (6.4)]{cv02}, which is still valid in the discrete-time case.  The first series is $\Theta$ by definition, the
second is equal to $\E_\pi [\theta_0]$ by \cite[Equation (6.6)]{cv02}.
It is obvious that $\E_\pi [\theta_0] \leq \Theta$.

If $\Theta < \infty$, then $\E_\pi [\theta_0] < \infty$, the deviation matrix exists and $K' = \Theta - \E_\pi [\theta_0] $ by~(\ref{e:one}).
If $\Theta = \infty$ and $\E_\pi [\theta_0] < \infty$, then $K'=\infty$ by (\ref{e:one}) again.  Finally, if $\E_\pi [\theta_0] = \infty$, then we have already seen in Remark \ref{r:deviation} that $K'=\infty$.
\end{proof}

\begin{remark}  \rm
First, let us deal with discrete-time birth and death processes. From \cite[Equation (6.6)]{cv02}, we see that $\Theta = \lim_{n\to \infty} E_n[\theta_0]$ and we may interpret Theorem~\ref{t:bdinfinite} as saying that, for Kemeny's constant to be finite, it is necessary (and sufficient) that having ventured to any state $n$, no matter how far from the origin, the process will reach state $0$ in bounded expected time.

In discrete-time, every transition from a state to one of its neighbours requires at least one unit of time, so that $\E_n[\theta_0] \geq n$ is unbounded.  This tells us that (\ref{e:sum}) diverges for all discrete-time birth-and-death processes.
\end{remark}

\begin{remark}   \rm
For the continuous-time birth and death process with birth rates $\lambda_n$ and death rates $\mu_n$, the right hand side of (\ref{e:G}) is the  `$D$ Series', see Anderson~\cite[Page 261]{ander91}, or Kijima~\cite [page 245]{kijim97}.  A recurrent continuous-time birth and death process with a finite $D$ Series is said to have an {\it entrance boundary at $\infty$}, a classification that goes back to Feller~\cite{felle59}. Bansaye, M\'el\'eard and Richard \cite{bamr16} described this as {\it instantaneously coming down from infinity}.

A number of authors have looked at consequences of the D-series being finite. For example, when there is an absorbing state at -1, this condition is equivalent to the existence of a unique quasistationary distribution (see van Doorn \cite[Theorem 3.2]{vand91}). The condition is also equivalent to {\it strong ergodicity} of the birth and death process in the sense that $\lim_{t\to \infty} \sup_i|p_{ij}(t) - \pi_j| = 0$, see \cite[Theorem 3.1]{zclh01}, \cite[Corollary 2.4]{zhan01} and \cite[Theorem 3.1]{mao02}.
\end{remark}

For a continuous-time birth and death process with finite $K'$, we can show that the change of limit and sum in  (\ref{e:indept}) can be justified, and so the physical interpretation given in Section \ref{s:finite} holds in this case as well. The details are given in the following lemma.
\begin{cor}
   \label{t:physInf}
Consider a continuous-time,  irreducible, positive recurrent, birth-and-death process on
$\{0, 1, \ldots\}$.  If $\Theta < \infty$, then
\[
K' = \lim_{t \rightarrow \infty}   (\sum_{j \in \s}\E_j[M_j(t)] - t).
\] 
\end{cor}

\begin{proof}
We fix $i=0$ and note that 
\begin{align*}
K' & = \sum_{j \in \s} \pi_j \E_0 [\theta_j] \\
  & = \sum_{j \in \s} \lim_{t \rightarrow \infty} (\E_j[M_j(t)] - E_0[M_j(t)])
\end{align*}
by Lemma \ref{t:continuous}.  We plan to use the Fatou-Lebesgue dominated convergence theorem to justify the change of limit and sum. To that end we construct a bound $m_{j}$ for $\E_j[M_j(t)] - E_0[M_j(t)]$ such that $\sum_j m_j < \infty$.  

Denote by $^{(0)}m_j(t)$ the expected sojourn time in $j$ during the interval $(0,t)$, starting from $j$, under taboo of state 0 and let $^{(0)}m_j( \infty )= \lim_{t \to \infty} (^{(0)}m_j(t))$.  Then
\begin{align*}
\E_j[M_j(t)] - E_0[M_j(t)] &= \, ^{(0)}m_j(t) + \int_0^t \E_0[M_j(t-u)]
                             \, \ud \P_j[\theta_0 \leq u]  -
                             \E_0[M_j(t)]
  \\
  &  \leq \, ^{(0)}m_j(t)  - (1-\P_0[\theta_0 \leq t]) \E_0[M_j(t)]
  \\
  & \leq \, ^{(0)}m_j(t)
  \\
  & \leq \, ^{(0)}m_j(\infty)
\end{align*}
It is a simple matter to show that $^{(0)}m_j(\infty) = \pi_j \sum_{1 \leq k \leq j} \frac{1}{\pi_k \mu_k}$ and that $\sum_{j \in \s} \,
^{(0)}m_j(\infty) = \Theta$, which is finite by assumption.  Thus,
\begin{align*}
K' & = \lim_{t \rightarrow \infty} \sum_{j \in \s}  (\E_j[M_j(t)] - E_0[M_j(t)])\\
& = \lim_{t \rightarrow \infty} \sum_{j \in \s}  (\E_j[M_j(t)] - t])
\end{align*}
by dominated convergence, and this concludes the proof.
\end{proof}

\begin{example}
   \label{e:mm1}
\rm
For the M/M/1 queue, $\lambda_n = \lambda$ and $\mu_n = \mu$,
independently of $n$.  The process is positive recurrent if and only
if the ratio $\rho = \lambda/\mu$ is strictly less than 1, and $\pi_n =
(1-\rho) \rho^n$.  Equation (\ref{e:G}) becomes
\begin{align*}
  \Theta &= \sum_{k \geq 0} \lambda^{-1}  \rho^{-k} \sum_{j \geq
           k+1} \rho^j  \\
  & = \sum_{k \geq 0} 1/(\mu - \lambda) = \infty
\end{align*}
so that Kemeny's constant is infinite.  However, it may be finite for
a process that we name the {\em sped-up}  M/M/1 queue: we
take an arbitrary sequence $\{\lambda_n\}$ and define $\mu_n = \rho
\lambda_{n-1}$, with $\rho < 1$.  Here, $\beta_n = \rho^n$ so that the
process is positive recurrent, $\pi_n = (1-\rho) \rho^n$ and 
\begin{align*}
  \Theta &= \sum_{k \geq 0} \lambda_k^{-1}  \rho^{-k} \sum_{j \geq
           k+1} \rho^j  \\
  & = \sum_{k \geq 0}   \lambda_k^{-1} \rho /(1-\rho)
\end{align*}
which converges if $\lambda_n \rightarrow \infty$ sufficiently fast.
In that case, $\mu_n$ tends to $\infty$ also.
\end{example}

The sped-up M/M/1 queue example illustrates that for Kemeny's
constant to be finite, transitions have to occur faster as the process
is further away from 0.  Actually, as we show in the next lemma, it is
necessary that transitions from $n$ to $n-1$ occur sufficiently fast,
transition rates from $n$ to $n+1$ being less critical.

\begin{lemma}
   \label{t:mun}
For $\Theta$ to be finite, it is necessary, but not sufficient, that
the series
$\sum_{j \geq 1} 1/\mu_j$ converges.
\end{lemma}
\begin{proof}
We rewrite (\ref{e:G}) as 
\begin{equation}
   \label{e:theta}
\Theta = \sum_{j \geq 1} f_j
\end{equation}
with
\begin{align}
  \nonumber
f_j & = \pi_j  \sum_{0 \leq k \leq j-1} (\lambda_k \pi_k)^{-1}  
 \\ \nonumber
  & = \pi_{j-1} \frac{\lambda_{j-1}}{\mu_j} (\sum_{0 \leq k \leq j-2}
    (\lambda_k \pi_k)^{-1}  + (\lambda_{j-1} \pi_{j-1})^{-1} )
 \\ \label{e:fj}
  & = (\lambda_{j-1} f_{j-1} +1) /\mu_j  \qquad \mbox{for $j \geq 1$,}
\end{align}
if we define $f_0=0$.

This shows that $f_j \geq 1/\mu_j$, so that the series (\ref{e:theta})
diverges if $\sum_{j \geq 1} 1/\mu_j$ diverges.
The proof that this is not a sufficient condition is given by
Example~\ref{e:counter} below.
\end{proof}

\begin{remark}    \rm
Lemma~\ref{t:mun} gives a different justification for the fact that (\ref{e:sum}) diverges for all discrete-time birth-and-death processes.: here, $\mu_n \leq 1-\lambda_n < 1$ by assumption, and the series $\sum_{j \geq 1} 1/\mu_j$ diverges.
\end{remark}

\begin{example}  
   \label{e:finiteK}
\rm
A direct consequence of Lemma \ref{t:mun} is that $K'$ is infinite for the
M/M/$\infty$ queue for which $\mu_n = n \mu$: the transition
rates from $n$ to $n-1$ are not large enough.  We may, however, use
Lemma~\ref{t:mun} to {\em design} processes for which Kemeny's constant is
finite.  To that end, we choose a sequence $\{f_j\}$ such that the series (\ref{e:theta}) converges, use (\ref{e:fj}) to define the sequence 
\begin{equation}
   \label{e:muj}
\mu_j =  (\lambda_{j-1} f_{j-1} +1) / f_j,
\end{equation}
and then find a sequence $\{\lambda_j\}$ such that the process is positive recurrent, that is,
such that $\sum_{n \geq 0} \beta_n$ converges, with $\beta_n$ defined
in (\ref{e:betan}).   Two such examples follow.  For the first,
\begin{itemize}
\item[] $f_1 = 0$, $f_j = 1/j^2$, for $j \geq 1$, 
\item[] $\mu_1=1$, $\mu_j = j^2(1+ 1/(j-1)^2)$, for $j \geq 2$, and
\item[] $\lambda_j = 1$ for all   $j$.
\end{itemize}
We easily see that 
\[
\beta_n = 1/ \prod_{2 \leq j \leq n}   j^2(1+1/(j-1)^2) < 1/(n!)^2.
\]
For the second, 
\begin{itemize}
\item[]
$f_j= \gamma^j$, with $\gamma < 1$,
\item[]
$\mu_j = \gamma^{-j} + \lambda_{j-1} \gamma^{-1}$, and
\item[]
$\{\lambda_j\}$ is arbitrary.
\end{itemize}
Here, 
\[
\beta_n = \prod_{0 \leq j \leq n-1} \lambda_j / \prod_{0 \leq j \leq
  n-1} (\lambda_j \gamma^{-1} + \gamma^{-(j+1)}) < \gamma^n.
\]
\end{example}

\begin{example}
   \label{e:counter}
\rm
This last example shows that Lemma \ref{t:mun} is not a necessary {\em
  and sufficient} condition.   Take
\begin{itemize}
\item[] 
$\mu_j = j^{1+\alpha}$, with $0 < \alpha < 1$, and
\item[] 
$\lambda_j = \mu_j$, for $j \geq 1$, $\lambda_0=1$.
\end{itemize}
With these parameters $\beta_n = 1/\mu_n$,  so that both $\sum_{n \geq
1} 1/\mu_n$ and $\sum_{n \geq 1} \beta_n$ converge.

By (\ref{e:fj}), we have $f_j \mu_j = 1 + f_{j-1} \mu_{j-1} = j$, so
that
\[
\sum_{j \geq 1} f_j = \sum_{j \geq 1} j/\mu_j = \sum_{j \geq 1} 1/j^\alpha
\]
diverges.
\end{example}

\section{Acknowledgment}
  \label{s:conclusion}

{This paper is an outgrowth of discussions that the second, third and fifth authors had during the International Workshop on Matrices and Statistics in
Funchal, in June, 2016. The authors would like to thank an anonymous referee for pointing out a couple of errors in an earlier version, as well as making a number of insightful comments.}

P.G. Taylor's research is supported by the Australian Research Council
(ARC) Laureate Fellowship FL130100039 and the ARC Centre of Excellence
for the Mathematical and Statistical Frontiers (ACEMS).

\appendix

\section{Proof that when $\cald$ exists, $K'$ is given by its trace}

Here we show that, whenever the deviation matrix $\cald$ exists, which occurs if and only if $\E_{\vpi}[\theta_j] < \infty$ for some state
$j \in \s$, then $K'$ is finite and equal to $\sum_{j \in \s} \cald_{jj}$ if this sum is finite, and infinite when $\sum_{j \in \s}  \cald_{jj}$ is infinite.

Since the deviation matrix exists, Coolen-Schrijner and van Doorn \cite[Equation (5.5)]{cv02} implies that
\begin{equation}
\label{e:dij}
\cald_{ij} = \pi_j\left(E_{\pi}[\theta_j] - E_i[\theta_j]\right)
\end{equation}
which gives us 
\begin{equation}
\label{e:pjeij}
\pi_j E_i[\theta_j] = \cald_{jj} - \cald_{ij},
\end{equation}
observing that $E_j[\theta_j]=0$. Summing (\ref{e:pjeij}) over $j$ we see that
\begin{align}
  \nonumber
\sum_{j \in \s} \pi_j \E_i[\theta_j] & = \sum_{j \in \s} \left(\cald_{jj} - \cald_{ij}\right)
 \nonumber \\
  & = \lim_{K\to \infty}\sum_{j \in \s_K} \left(\cald_{jj} - \cald_{ij}\right)
 \nonumber \\
& =  \lim_{K\to \infty}\left(\sum_{j \in \s_K} \cald_{jj} - \sum_{j \in \s_K}\cald_{ij}\right),
 \label{e:Kinfinite}
 \end{align}
where $\{\s_K\}$ is a monotone sequence of finite subsets converging to $\s$. 

By \cite[Theorem 5.2]{cv02}, $\sum_{j \in \s}\cald_{ij} = 0$ for all $i \in \s$. Therefore, for any $\epsilon > 0$, there exists $K_0(i)$ such that 
\[
-\epsilon/2 \leq \sum_{j \in \s_K}\cald_{ij} \leq \epsilon/2
\]
for all $K \geq K_0(i)$.

If $\sum_{j \in \s} \cald_{jj} = \infty$ then, for any $M>0$, there exists $K_1$ such that $\sum_{j \in \s_K} \cald_{jj} \geq M$ for all $K > K_1$. It follows that, if $K > \max(K_0(i),K_1)$ then 
\[
\sum_{j \in \s_K} \cald_{jj} - \sum_{j \in \s_K}\cald_{ij} \geq M - \epsilon/2
\]
and so 
\[
\sum_{j \in \s} \pi_j \E_i[\theta_j] = \sum_{j \in \s} \left(\cald_{jj} - \cald_{ij}\right) = \infty
\]
independently of $i$.

On the other hand, if $\sum_{j \in \s} \cald_{jj} = L < \infty$ then, for any $\epsilon>0$, there exists $K_1$ such that $L - \epsilon/2 \leq \sum_{j \in \s_K} \cald_{jj} \leq L + \epsilon/2$ for all $K > K_1$. Then, again for fixed $i$, taking $K > \max(K_0(i),K_1)$,
\[
 L - \epsilon \leq \sum_{j \in \s_K} \cald_{jj} - \sum_{j \in \s_K} \cald_{ij}\leq L + \epsilon
\]
and so
\[
\sum_{j \in \s} \pi_j \E_i[\theta_j] = \sum_{j \in \s} \cald_{jj} = L.
\]
This argument holds for all $i$, and we see that $\sum_{j \in \s} \pi_j \E_i[\theta_j]$ is independent of~$i$.

\end{document}